\documentclass[11pt,reqno]{amsart}
  
\usepackage{latexsym}
\usepackage{graphicx}
\usepackage{amssymb}
 \graphicspath{{Figures/}}
 \usepackage[show]{ed}

\newcommand{\red}[1]{{#1}}

\renewcommand{\Re}{{\operatorname{Re}\,}}

\renewcommand{\epsilon}{\varepsilon}

\newcommand{\R}{{\mathbb R}}

\newcommand{\Z}{{\mathbb Z}}

\newcommand{\lan}{\left\langle}
\newcommand{\ran}{\right\rangle}
\newcommand{\mc}[1]{\mathcal{#1}}
\newcommand{\e}{\epsilon}

\newcommand{\re}{\mathbb{R}}

\newcommand{\supp}{{\operatorname{supp\,}}}
\renewcommand{\phi}{\varphi}

\newcommand{\ep}{\varepsilon}

\newtheorem{theo}{{\sc Theorem}}

\newtheorem{cor}{{\sc Corollary}}[section]

\newtheorem{lem}[cor]{{\sc Lemma}}

\numberwithin{equation}{section}

\newenvironment{rem}[1][]{\refstepcounter{cor}{\medskip\noindent{\it Remark~\thecor :\/#1} }}{\medskip}

\newcommand{\spa}{\operatorname{span}}
\newtheorem{defn}[cor]{{\sc Definition}}

\usepackage{hyperref}
\newcommand{\admissible}{{admissible }}

\title[Eigenfunction bounds]{Eigenfunction scarring and improvements in $L^{\infty}$ bounds}

\author{Jeffrey Galkowski}
\address{Department of Mathematics and Statistics, McGill University, Montr\'eal, QC, Canada}
\email{jeffrey.galkowski@mcgill.ca }
\author{John A. Toth}
\address{Department of Mathematics and Statistics, McGill University, Montr\'eal, QC, Canada}
\email{jtoth@math.mcgill.ca} 
\date{}

\begin{document}

\begin{abstract}
We study the relationship between $L^\infty$ growth of eigenfunctions and their $L^2$ concentration as measured by defect measures. In particular, we show that scarring in the sense of concentration of defect measure on certain submanifolds is incompatible with maximal $L^\infty$ growth. In addition, we show that a defect measure which is too diffuse, such as the Liouville measure, is also incompatible with maximal eigenfunction growth. 
\end{abstract}
\maketitle
\section{Introduction}

Let $(M,g)$ be a $C^{\infty}$ compact manifold of dimension $n$ without boundary. Consider the eigenfunctions 
\begin{equation}
\label{e:eigenfunction}
(-\Delta_g-\lambda_j^2)u_{\lambda_j}=0,\quad \|u_{\lambda_j}\|_{L^2}=1
\end{equation}
as $\lambda_j\to\infty$. It is well known \cite{Ava,Lev,Ho68} (see also \cite[Chapter 7]{EZB}) that solutions to~\eqref{e:eigenfunction} satisfy
\begin{equation}
\label{e:L infinity}
\|u_{\lambda_j}\|_{L^\infty(M)}\leq C\lambda_j^{\frac{n-1}{2}}
\end{equation}
and that this bound is saturated e.g. on the sphere. It is natural to consider the situations which produce sharp examples for (\ref{e:L infinity}).  In many cases, one expects polynomial improvements  to (\ref{e:L infinity}), but rigorous results along these lines are few and far between \cite{I-S}. \red{In the case of negatively curved manifolds, $\log$ improvements can be obtained~\cite{Berard77}. However,}  at present, under general dynamical assumptions, known results involve $o$-improvements to (\ref{e:L infinity}) \cite{TZ02,SoggeZelditch,TZ03,SoggeTothZelditch,SZ16I,SZ16II}.
These papers all study the  connections between the growth of $L^\infty$ norms of eigenfunctions and the global geometry of the manifold $(M,g)$. 
In this note, we examine the relationship between $L^\infty$ growth and $L^2$ concentration of eigenfunctions. 
We measure $L^2$ concentration using the concept of a {\em{defect measure}} - a sequence $\{u_{\lambda_j}\}$ has defect measure $\mu$ if for any $a\in S_{\text{hom}}^0(T^*M\setminus \{0\})$, 
\begin{equation} \label{defect}
\lan a(x,D) u_{\lambda_j},u_{\lambda_j}\ran \to \int_{S^*M}  a(x,\xi)d\mu. \end{equation}

By an elementary compactness/diagonalization argument it follows that any sequence of eigenfunctions $u_{\lambda_j}$ solving~\eqref{e:eigenfunction} possesses a further subsequence that has  a defect measure in the sense of~\eqref{defect} (\cite[Chapter 5]{EZB},\cite{Ger91}).  Moreover, a standard commutator argument shows that if $\{ u_{\lambda_j} \}$ is any sequence of $L^2$-normalized Laplace eigenfunctions, the associated defect measure
  $\mu$ is invariant under the geodesic flow; that is, if $G_t: S^*M  \to S^*M$ is the geodesic flow \red{ i.e. the hamiltonian flow of $p=\frac{1}{2}|\xi|^2_g$},
$ (G_t)_*\mu =\mu,\,\forall t \in \R$. 

\begin{defn}
We say that an eigenfunction subsequence is \emph{strongly scarring} provided that \red{for any defect measure $\mu$ associated to the sequence,} $\supp \mu$ is a finite union of periodic geodesics.  
\end{defn}
\begin{theo}
\label{mainthmLaplace}
Let $\{u_{\lambda_j}\}$ be a strongly scarring sequence of solutions to~\eqref{e:eigenfunction}. Then
$$\|u_{\lambda_j}\|_{L^\infty}=o(\lambda_j^{\frac{n-1}{2}}).$$
\end{theo}
We also have improved $L^\infty$ bounds when eigenfunctions are \emph{quantum ergodic}, that is, their defect measure is the Liouville measure on $S^*M$, $\mu_L$ \red{(see e.g.~\cite{SniQE, CdVQE, ZeQE} for the standard quantum ergodicity theorem)}. 
\begin{theo} \label{QE}
Let $\{u_{\lambda_j}\}$ be a quantum ergodic sequence of solutions to~\eqref{e:eigenfunction}. Then
$$\|u_{\lambda_j}\|_{L^\infty}=o(\lambda_j^{\frac{n-1}{2}}).$$ 
\end{theo}

Theorems~\ref{mainthmLaplace} and~\ref{QE} are corollaries of our next theorem where we relax the assumptions on $\mu$ and make the following definitions. Define the \emph{time $T$ flow out} by
$$\Lambda_{x,T}:=\bigcup_{t=-T}^TG_t(S^*_xM).$$

\begin{defn} 
Let ${\mathcal H}^{n}$ be $n$-dimensional Hausdorff measure on $S^*M$ induced by the Sasaki metric on $T^*M$ (see for example \cite[Chapter 9]{BlairSasaki} for a treatment of the Sasaki metric).
 We say that the subsequence $u_{\lambda_j}; j=1,2,...$ is \emph{\admissible at $x$} if \red{for any defect measure $\mu$ associated to the sequence} \red{ there exists $T>0$ such that }
\begin{equation} \label{admissible}
{\mathcal H}^n (  \,\supp \,\, \mu |_{\Lambda_{x\red{,T}}} \, ) = 0. \end{equation}
We say that the subsequence is \emph{\admissible}if it is admissible at $x$ for every $x\in M$.
\end{defn}

We note that in (\ref{admissible})  $\mu |_{\Lambda_{x\red{, T}}}$ denotes the defect measure restricted to the flow out $\Lambda_{x\red{, T}};$ for any $A$ that is $\mu$-measurable,
$$\mu|_{\Lambda_{x\red{, T}}}(A):= \mu ( A \cap \Lambda_{x\red{, T}}).$$\

\begin{theo}\label{thmweakscarring}
\label{mainthm}
Let $\{ u_{\lambda_j} \}$ be a sequence of $L^2$-normalized Laplace eigenfunctions that is admissible in the sense of~\eqref{admissible}. 
Then 
$$\|u_{\lambda_j}\|_{L^\infty}=o(\lambda_j^{\frac{n-1}{2}}).$$
\end{theo}

\begin{rem}
We choose to use the Sasaki metric to define $\mathcal{H}^n$ for concreteness, but this is not important and we could replace the Sasaki metric by any other metric on $S^*M$. 
\end{rem}


Theorem~\ref{mainthm} can be interpreted as saying that eigenfunctions which strongly scar are too concentrated to have maximal $L^\infty$ growth, while diffuse eigenfunctions are too spread out to have maximal growth. However, the reason the adimissiblity assumption is satisfied differs in these cases. In the diffuse case (see Theorem~\ref{QE}), one has $ \mu |_{\Lambda_{x,\red{T}}} = 0,$ so that the admissibility assumption is trivially verified. In the case where the eigenfunctions strongly scar (see Theorem~\ref{mainthmLaplace}), $\mu|_{\Lambda_{x,\red{T}}} \neq 0$ but the Hausdorff dimension of $\supp \, \mu|_{\Lambda_{x,\red{T}}}$ is $< n;$ so again, (\ref{admissible}) is satisfied. The zonal harmonics on the sphere $S^2$\red{, which saturate the $L^\infty$ bound~\eqref{e:L infinity},}  lie precisely between being diffuse and strongly scarring (see section~\ref{sec:zone}). 
 
 Observe that the condition $\mu$ is diffuse is much more general than $\mu=\mu_{L}$. Jakobson--Zelditch~\cite{JaZe} show that any invariant measure on $S^*S^n$ where $S^n$ is the round sphere can be obtained as a defect measure for a sequence of eigenfunctions and in particular many non-Liouville but diffuse measures occur.\\

\begin{rem}
We note that the results here hold for any $o(\lambda)$ quasimode of $(-\Delta_g-\lambda^2)$ that is compactly microlocalized in frequency (see~\cite{GLinfty}).
\end{rem}

\subsection{Relation with previous results}
Theorem~\ref{QE} is related to \cite[Theorem 3]{SoggeTothZelditch}, where the $o(h^{\frac{1-n}{2}} )$ sup bound is proved for all Laplace eigenfunctions on a $C^{\omega}$ surface with ergodic geodesic flow. However, in Theorem~\ref{QE}, we make no analyticity or dynamical assumptions on $(M,g)$ whatsoever, only an assumption on the particular defect measure associated with the eigenfunction sequence. Recently, Hezari \cite{Hez16} \red{and Sogge \cite{SoggeLocal}} gave independent proofs of Theorem~\ref{QE}.

\red{One consequence of the work of Sogge  is  the relation between $L^p$ norms for eigenfunctions and the push forward of defect measures to the base manifold $M$.  In particular, he shows~\cite[(3.3)]{SoggeLocal} that
\begin{equation} \label{soggebound}
 \| u_\lambda \|_{L^{\infty}(M)} \leq C \lambda^{\frac{n-1}{2}} \sup_{x \in M} \delta^{-1/2} \| u_{\lambda} \|_{L^2(B_{\delta}(x))} 
 \end{equation}
when $\lambda^{-1} \leq \delta \leq inj(M,g)$ and $\lambda \geq 1.$ 
We note that when $u_{\lambda}$ are quantum ergodic, $\| u_{\lambda} \|_{L^2(B_\delta(x)} \approx \delta^{\frac{n}{2}}$ and so the $o(\lambda^{\frac{n-1}{2}})$-bound in Theorem \ref{QE} follows from (\ref{soggebound}) as well (see also Corollary 1.2 in \cite{SoggeLocal}).}

\red{ 
However, neither the scarring result in Theorem \ref{mainthmLaplace} nor the more general bound in Theorem \ref{thmweakscarring} follow from~\eqref{soggebound}.  To compare and contrast with (\ref{soggebound}), we observe that \eqref{soggebound} implies for any $\delta>0$ independent of $\lambda$,
$$
\limsup_{\lambda \to \infty }\lambda^{\frac{1-n}{2}} \| u_\lambda \|_{L^{\infty}(M)}\leq C \sup_{x \in M} \delta^{-1/2} \, \Big( \mu( S^*B_{\delta}(x) ) \Big)^{\frac{1}{2}}.
$$
 Our main estimate in (\ref{e:related})  says that for any $x(\lambda )$ with $d(x(\lambda ),x)=o(1)$,
\begin{equation} \label{mainbound}
\limsup_{\lambda \to \infty }\lambda^{\frac{1-n}{2}} |u_{\lambda}(x(\lambda))|\leq  C'_{\delta}\Big(\mathcal{H}^n(\supp \mu|_{A_x(\delta/2,3\delta)})\Big)^{1/2}
\end{equation}
where for $\delta_2 > \delta_1$ the $A_x(\delta_1,\delta_2) = \Lambda_{x,\delta_2} \setminus \Lambda_{x,\delta_1}.$ This microlocalized bound allows us to deal with the more general scarring-type cases as well. In particular, the key differences are that we have replaced $S^*B_\delta(x)$ by $A_x(\delta/2,2\delta)\subset \Lambda_x$ and the defect measure by Hausdorff $n$ measure. We note however that unlike (\ref{soggebound}),  $\delta >0$ can be arbitrarily small but is fixed independent of $\lambda$ in (\ref{mainbound}).}

In \cite{SoggeZelditch}, Sogge--Zelditch prove that any manifold on which~\eqref{e:L infinity} is sharp must have a self focal point. That is, a point $x$ such that $|\mc{L}_x|>0$ where 
$$\mc{L}_x:=\{\xi\in S^*_xM\mid \text{ there exists }T\text{ such that }\exp_xT\xi=x\}$$
and $|\cdot|$ denotes the normalized surface measure on the sphere. Subsequently, in \cite{SoggeTothZelditch} the authors showed that one can replace $\mc{L}_x$ by the set of recurrent directions $\mc{R}_x \subset \mc{L}_x$ and the assumption $|\mc{R}_x|>0$ for some $x \in M$ is necessary to saturate the maximal bound in (\ref{e:L infinity}). \red{Here,
$$\mc{R}_x:=\Big\{\xi \in S^*_xM\mid \xi\in \Big(\bigcap_{T>0} \overline{\bigcup_{t\geq T}G_t(x,\xi)\cap S^*_xM}\Big)\bigcap \Big(\bigcap_{T>0} \overline{\bigcup_{t\leq -T}G_t(x,\xi)\cap S^*_xM}\Big)\Big\}.$$
} The example of the triaxial ellipsoid with $x$ equal to an umbilic point shows that latter assumption is weaker than the former. Indeed, in such a case $|\mc{L}_x|=1$ whereas $|\mc{R}_x|=0.$  Most recently, in \cite{SZ16I,SZ16II}, it was proved that for real-analytic surfaces, the maximal $L^\infty$ bound can only achieved if there exists a periodic point $x \in M$ for the geodesic flow, \red{i.e. a point $(x,\xi)$ so that} all geodesics starting at $(x,\xi) \in S^*M$ close up smoothly after some finite time $T>0.$

Together with our analysis, the results of \cite{SoggeTothZelditch} imply that any sequence of eigenfunctions, $\{u_\lambda\}$ having maximal $L^\infty$ growth \red{near $x$} and defect measure $\mu$ must have \red{$\mu(\Lambda_{x,T})>0$ for all $T>0$ and $|\mathcal{R}_x|>0$.} \red{In particular, the results of \cite{SoggeTothZelditch} show that $u_\lambda$ can only have maximal $L^\infty$ growth near a point with a positive measure set of recurrent points and Theorem~\ref{mainthm} shows that a point with maximal $L^\infty$ growth must have $\mu(\Lambda_{x,T})>0$.}
As far as the authors are aware, the results in \cite{SoggeTothZelditch} and in \cite{SZ16I, SZ16II} do not give additional information about $\mu$. 

On the other hand, under an additional regularity assumption on the measure $\mu$, Theorem~\ref{mainthm} \red{can be used to show that when $u_\lambda$ has maximal growth near $x$,} $\mu|_{\Lambda_{x.\red{T}}}$ is not mutually singular with respect to $\mc{H}^n$. \red{Since the measure for a zonal harmonic is a smooth multiple of $\mc{H}^n$  (see Section~\ref{sec:zone})}, this implies that the measure $\mu$ resembles the defect measure of a zonal harmonic . In~\cite{GLinfty}, the first author removes the necessity for any additional regularity assumption and gives a full characterization of defect measures for eigenfunctions with maximal $L^\infty$ growth\red{, in particular proving that if $u_\lambda$ has maximal growth near $x$ and defect measure $\mu$, then $\mu|_{\Lambda_{x,T}}$ is not mutually singular with respect to $\mc{H}^n$}. Finally, we note that unlike \cite{SoggeZelditch,SoggeTothZelditch,SZ16I,SZ16II}, the analysis here is entirely local.\\

\noindent {\sc Acknowledgemnts.} The authors would like to thank the anonymous referees for their detailed reading and many helpful comments. J.G. is grateful to the National Science Foundation for support under the Mathematical Sciences Postdoctoral Research Fellowship  DMS-1502661.  The research of J.T. was partially supported by NSERC Discovery Grant \# OGP0170280 and an FRQNT Team Grant. J.T. was also supported by the French National Research Agency project Gerasic-ANR-
13-BS01-0007-0.


\section{A local version of~\ref{thmweakscarring}}

In the following, we will freely use semiclassical pseudodifferential calculus where the semiclassical parameter is $h$ with $h^{-1} = \lambda \in \text{Spec} \, \sqrt{-\Delta_g}.$ \red{We write $r(x,y):M\times M\to \mathbb{R}$ for the Riemannian distance from $x$ to $y$ and write $B(x,\delta)$ for the geodesic ball of radius $\delta$ around $x$.} We start with a local result:
 \begin{theo} \label{thm:local}
Let $\{ u_{h} \}$ be sequence of Laplace eigenfunctions that is \admissible  at $x$.
Then for any $\red{\delta}(h)=o(1)$, 
$$\|u_h\|_{L^\infty(B(x,\red{\delta}(h))}=o(h^{\frac{1-n}{2}}).$$
\end{theo}
Theorem~\ref{mainthm} is an easy consequence of Theorem~\ref{thm:local}.
\begin{proof}[Proof that Theorem~\ref{thm:local} implies Theorem~\ref{mainthm}]
Suppose that $u$ is \admissible  and 
$$\limsup_{h\to 0} h^{\frac{n-1}{2}}\|u_h\|_{L^\infty}\neq 0.$$
Then, there exist $c>0$, $h_k\to 0$, $x_{h_k}$ so that 
$$|u_{h_k}(x_{h_k})|\geq ch_k^{-\frac{n-1}{2}}.$$
Since $M$ is compact, by taking a subsequence, we may assume $x_{h_k}\to x$. But then $\red{r}(x,x_{h_k})=o(1)$ and since $u$ is \admissible at $x$, Theorem~\ref{thm:local} implies
$$\limsup_{k\to \infty} h_k^{\frac{n-1}{2}}|u_{h_k}(x_{h_k})|=0.$$
\end{proof}



\section{Proof of Theorem~\ref{thm:local} }
\label{sec:easyDiffuse}

In view of the above, it suffices to prove the local result: Theorem~\ref{thm:local}.

\begin{proof}

Fix $T>3\delta >0$ and let $\rho \in \mc{S}(\R)$ with $\rho(0)=1$ and 
$\supp\hat{\rho} \subset (\delta, 2 \delta).$
Let 
$$S^*M(\e):= \{ (x,\xi); | |\xi|_{x} - 1| \leq \e \}$$
 and  $\chi(x,\xi) \in C^{\infty}_{0}(T^*M)$ be a cutoff near the cosphere $S^*M$ with
$ \chi(x,\xi) =1$ for $(x,\xi) \in S^*M(\ep)$ and $\chi(x,\xi) = 0$ when $(x,\xi) \in T^*M \setminus S^*M(2\e).$ Let  $\chi(x,hD) \in Op_h(C_0^\infty(T^*M))$ be the corresponding $h$-pseudodifferential cutoff. Also, in the following, we will use the notation
$$ \Gamma_x := \supp \,\, \mu|_{\Lambda_{x\red{, T}}}$$
to denote the support of the restricted defect measure corresponding to the eigenfunction sequence $\{ u_{h_j} \}$ in Theorem~\ref{mainthm}.

Then,  we have
\begin{align} \label{QE0}
u_h &= \rho\Big( \red{\frac{1}{2h}} [ -h^2 \Delta -1]\Big) u_h = \int_{\R} \hat{\rho}(t)  e^{i \red{\frac{t}{2}}[-h^2\Delta-1]/h} \chi(y,hD_y) u_h \, dt + O_{\ep}(h^{\infty}). \end{align}

\subsection{Microlocalization to the flow out $\Lambda_x$}

Set 
$$V(t,x,y,h):=  \Big( \hat{\rho}(t)  e^{i\red{\frac{t}{2}} [-h^2\Delta-1]/h} \chi(y,hD_y) \Big) (t,x,y).$$
Then, \red{by Egorov's Theorem~\cite[Theorem 11.1]{EZB}}
\begin{equation}
\label{e:octopoda}
WF_h'( V(t,\cdot, \cdot,h)) \subset \{ (x,\xi,y,\eta); (x,\xi) = G_t(y,\eta),  \,  | |\xi|_x -1 | \leq 2\ep \, \}.
\end{equation}
\red{(see e.g.~\cite[Definition E.37]{ZwScat} for a definition of $WF_h'$)}.
\\
Let $b_{x,\e}\red{(x,hD)} \in \red{Op_h}(C_0^\infty(T^*M))$ be a family of $h$-pseudodifferential cutoffs with \red{principal} symbols 
$$b_{x,\e}\in C^\infty_{0}(\{(y,\eta)\mid (y,\eta)=G_t(x_0,\xi)\, \text{for some} \, (x_0,\xi)\in S_{x_0}^*M(3\ep)\,\text{with}\, r(x,x_0)<2\e,  \frac{\delta}{2} < t < 3\delta\},$$ 
with 
$$b_{x,\e}\equiv 1\text{ on }\{(y,\eta)\mid (y,\eta)=G_t(x_0,\xi)\, \text{for some} \, (x_0,\xi)\in S_{x_0}^*M(2\ep)\,\text{with}\, r(x,x_0)<\e,  \delta < t < 2\delta\}.$$

By \red{the definition of $WF_h'$} together with\red{~\eqref{QE0} and}~\eqref{e:octopoda}, it follows that for $r(x(h),x)=o(1)$,
\begin{equation} \label{microlocalize-t}
u_h(x(h)) = \int_{M}  \bar{V}(x(h),y,h)  \, b_{x,\e}(y,hD_y) u_h(y) dy  + O_{\ep}(h^{\infty}), \end{equation}
where,
$$\bar{V}(x(h),y,h) := \int_{\R} \hat{\rho}(t) \big( e^{i\red{\frac{t}{2}} [-h^2\Delta-1]/h} \chi(y,hD_y) \big) (t,x(h),y) \, dt.$$\

By a standard stationary phase argument,
\begin{equation} \label{wkb-t}
\begin{gathered} 
\bar{V}(x,y,h) = h^{\frac{1-n}{2}} e^{ -i r(x,y))/h}   a(x,y,h) \,  \hat{\rho}(r(x,y)) + O_{\ep}(h^{\infty}), \end{gathered}\end{equation}\\
where $ a(x,y,h) \in S^{0}(1)$. \\

\red{To see this, observe that by~\cite[Theorem 10.4]{EZB}
$$
\bar{V}(x,y,h)=(2\pi h)^{-n}\int e^{i\varphi(t,x,y,\eta) /h}\,  \alpha(t,x,y,\eta,h)\hat{\rho}(t)d\eta dt+O(h^\infty)
$$
where $b \in C_c^\infty$ and $\varphi$ solves
\begin{equation}
\label{e:phase}
 \partial_t\varphi =\frac{1}{2}(|\partial_x\varphi|^2_{g(x)}-1),\qquad \varphi(0,x,y,\eta)=\langle x-y,\eta\rangle
 \end{equation}
 In particular, for all $(t,x,y,\eta)$, $\exp(tH_{\frac{1}{2}|\xi|_g^2})(\partial_\eta \varphi+y,\eta)=(x,\partial_{x}\varphi)$. The phase function
 $$\phi(t,x,y,\eta) = \langle \exp_{y}^{-1}(x),\eta \rangle + \frac{t}{2} (|\eta|^2_{y} -1)$$ 
 satisfies~\eqref{e:phase}.}

\red{  We next perform stationary phase in $(t,\eta)$.
 First, observe that the phase is stationary at 
 $$
 \exp(tH_{\frac{1}{2}|\xi|_g^2})(y,\eta)=(x,\partial_{x}\varphi),\qquad |\partial_x\varphi|_{g(x)}=1.
 $$
 In particular, $t=r(x,y)$ and the geodesic through $(y,\eta)$ passes through $x$. Since $ \supp \hat{\rho}\subset (\delta,2\delta)$, by performing non-stationary phase, we may assume $t\in (\delta,2\delta)$ and hence $\delta<r(x,y)<2\delta$. Then, we observe that $\partial^2_{(t,\eta)}\varphi$ is non-degenerate for $t\in (\delta,2\delta).$ The solutions $(t_c,\eta_c)$ of the critical point equations $\partial_t \phi =0$ and $\partial_{\eta} \phi = 0$ are given by 
 $$t_c =   |\exp_y^{-1}(x)| = r(x,y), \quad \,\,\, \eta_c = - \frac{\exp_y^{-1}(x)}{r(x,y)}.$$
 Consequently, \eqref{wkb-t} follows from an application of stationary phase. (see also~\cite[\red{Lemma 5.1.3}]{SoggeBook} or \cite[Theorem 4]{BGT}). }

Then, in view of (\ref{wkb-t}) and (\ref{microlocalize-t}),

\begin{equation} \label{QE1}
\begin{gathered}
u_h(x(h)) = v_h(x(h))+O_{\ep}(h^\infty)\\ \\
v_h(x(h))=h^{\frac{1-n}{2}} \int_{ \frac{\delta}{2} < \red{r(x,y)}< 2\delta} e^{ -ir(x(h),y)/h}  a(x(h),y,h) \hat{\rho}(r(x(h),y)) \, b_{x,\e}(y,hD_y) u_h(y) dy .\end{gathered} \end{equation}

\red{
Now, note that for any $\psi\in C_0^\infty(M)$, 
\begin{equation}
\label{e:decompU}
v_h(x(h))=I_{1}(x(h),h) + I_{2}(x(h),h)
\end{equation}
where
\begin{align} 
 I_1&:= (2\pi h)^{\frac{1-n}{2}}  \int_{ \delta/2 < \red{r(x,y)} < 2\delta} e^{ -ir(x(h),y)/h}  a(x(h),y,h) \hat{\rho}(r(x(h),y)) \, \psi(y) \, \,\red{(b_{x,\e}(y,hD_y)u_h)} \, dy   \nonumber \\
 \nonumber\\
I_2&:= (2\pi h)^{\frac{1-n}{2}} \int_{ \delta/2 < \red{r(x,y)} < 2\delta}  e^{ -ir(x(h),y)/h} a (x(h),y,h) \hat{\rho}(r(x(h),y)) \, (1- \psi(y) )\, \, \red{(b_{x,\e}(y,hD_y) u_h)} \, dy.\nonumber
 \end{align}\\
 Therefore,
by Cauchy-Schwarz applied to $I_1$ and $I_2$, 
$$|h^{\frac{n-1}{2}}v_h(x(h))|\leq C_\delta (\|\psi\|_{L^2}\|b_{x,\e}(y,hD_y)u_h(y)\|_{L^2}+\|(1-\psi(y))[b_{x,\e}(y,hD_y)u_h]\|_{L^2}).$$
Hence letting $h\to 0$ then $\e\to 0$, and using that (see for example~\cite[Theorem 5.1]{EZB})
$$\|b_{x,\e}(y,hD_y)u_h(y)\|_{L^2}\leq (\sup|b_{x,\e}|+o_\e(1))\|u_h\|_{L^2}$$
we have
\begin{equation}
\label{e:imp}
\limsup_{h\to 0}h^{\frac{n-1}{2}}|u_h(x(h))|\leq C_\delta \Big(\|\psi\|_{L^2}+\Big(\int_{\Lambda_{x,3\delta}\setminus \Lambda_{x,\delta/2}}(1-\psi(y))^2d\mu\Big)^{\frac{1}{2}}\Big)
\end{equation}
}

\subsection{Further microlocalization along $\supp \mu |_{\Lambda_x}$}

Let ${\mathcal H}^{n}$ be the $n$-dimensional Hausdorff measure on the flow out $\Lambda_{x}.$
By assumption, ${\mathcal H}^{n}( \supp \, \mu|_{\Lambda_x} ) =0.$ In view of the microlocalization above, we are only interested in the annular subset 
$$A_x(\delta/2,3\delta):= 
\Lambda_{x,3\delta} \setminus \Lambda_{x,\delta/2}.$$
\red{Since $\mc{H}^n$ is Radon,} for any $\e_1>0$, there exist $n$-dimensional balls $B(r_j) \subset A_x(\delta/4, 4\delta); j=1,2,...$ with radii $r_j>0, j=1,2,...$ 
such that 
$$ \supp \,  \mu|_{\red{A_x(\delta/2,3\delta)}} \subset \bigcup_{j=1}^\infty B(r_j), \quad {\mathcal H}^n\Big( \bigcup_{j=1}^\infty B(r_j) \Big) <\red{\mc{H}^n(\supp \, \mu|_{A_x(\delta/2,3\delta)})+} \ep_1.$$

Note that for $\delta>0$ small enough, the canonical projection $\pi: T^*M \to M$ restricts to a diffeomorphism
$$\pi:  A_x(\delta/4, 4\delta) \to \{ y \in M; \delta/4 < r(x,y) < 4\delta \}. $$

Consider the closed set 
$$K = \pi ( \supp \, \mu|_{  A_x(\red{\delta/2, 3\delta}) } ) \subset M \ $$\
\noindent with open covering
\begin{equation}
\label{open estimate}G:= \pi \Big( \bigcup_{j=1}^\infty B(r_j)\Big),\quad \text{ satisfying} \quad
{\mathcal H}^{n}(G) = \red{\mathcal{H}^n(K)}+O(\ep_1). \end{equation}

By the $C^{\infty}$ Urysohn lemma, there exists $\chi_{\Gamma_x} \in C^{\infty}_0(M; [0,1])$ with 
\begin{equation}
\label{e:chimp}
\chi_{\Gamma_x} |_K = 1,\qquad \supp \chi_{\Gamma_x} \subset G.
\end{equation}
(Note that $\chi_{\Gamma_x}$ depends on $\e_1$, but we suppress this dependence to simplify notation.) 
\red{
We now apply~\eqref{e:imp} with $\psi=\chi_{\Gamma_x}$. First, observe that by~\eqref{open estimate} and~\eqref{e:chimp}
\begin{equation}
\label{e:vol}
\|\chi_{\Gamma_x}\|_{L^2}\leq \big(\mc{H}^n(G)\big)^{1/2}\leq \big(\mc{H}^n(K)\big)^{1/2}+O(\e_1^{1/2}).
\end{equation}
}

Next, by construction, for all $\e_1>0$, 
$$(1-\chi_{\Gamma_x})(y) = 0, \quad \forall y \in \pi (  \, \supp \, \mu |_{ \Lambda_{x,4\delta} \setminus \Lambda_{x,\delta/4} } )$$
and hence 
$$\int_{\Lambda_{x,3\delta}\setminus \Lambda_{x,\delta/2}} (1-\chi_{\Gamma_x})^2d\mu=0.$$
Using this together with~\eqref{e:vol} in~\eqref{e:imp} and sending $\e_1\to 0$ gives 
\begin{equation}
\label{e:related}
\limsup_{h\to 0}h^{\frac{n-1}{2}}|u_{h}(x(h))|\leq C_\delta\big(\mc{H}^n(\pi( \supp \, \mu|_{  A_x(\delta/2, 3\delta) }))\big)^{\frac{1}{2}}\leq C'_{\delta}\Big(\mathcal{H}^n(\supp \mu|_{A_x(\delta/2,3\delta)})\Big)^{1/2}.
\end{equation}
where the last inequality follows from the fact that $\pi|_{A(\delta/2,3\delta)}$ is a diffeomorphism. Finally, since $u_h$ is admissible at $x$, 
$$\mathcal{H}^n(\supp \mu|_{A_x(\delta/2,3\delta)})=0$$
finishing the proof.\\

\begin{rem}
For $r(x(h),x)=o(1)$, the estimate
\begin{equation*}
\limsup_{h\to 0}h^{\frac{n-1}{2}}|u_{h}(x(h))|\leq  C'_{\delta}\Big(\mathcal{H}^n(\supp \mu|_{A_x(\delta/2,3\delta)})\Big)^{1/2}
\end{equation*}
in~\eqref{e:related} holds for any sequence of eigenfunctions with defect measure $\mu$. It gives a quantitative estimate relating the behavior of the defect measure to $L^\infty$ norms of eigenfunctions.
This estimate can also be obtained as a consequence of~\cite[Theorem 2]{GLinfty} by replacing the absolutely continuous part of $\mu$ with $1_{\supp \mu|_{\Lambda_x}}d\mc{H}^n.$
\end{rem}
\end{proof}

\section{The example of zonal harmonics}
\label{sec:zone}

Let $(S^2,g_{can})$ be the round sphere and $(r,\theta)$ be polar variables centered at the north pole $p = (0,0,1) \in \R^3.$ The geodesic flow is a completely integrable system with Hamiltonian
\begin{equation} \label{polar}
H = |\xi|_g^2 = \xi_r^2 + (\sin r)^{-2} \xi_{\theta}^2, \quad r \in (0,\pi) \end{equation}
and Claurault integral $p = \xi_{\theta}$
satisfying $\{ H, p \} = 0.$ The associated moment mapping is
${\mathcal P} = (H,p): T^*S^2 \to \R^2$ and the connected components of the level sets are, by the Liouville-Arnold Theorem, Lagrangian tori $\Lambda_{c}$ indexed by the values of the moment map $(1,c)  \in {\mathcal P}(T^*S^2).$

The associated quantum integrable system is given by 
the Laplacian $\Delta_g$ and the rotation operator $h D_{\theta}.$
The corresponding $L^2$-normalized joint eigenfunctions are the standard spherical harmonics $Y^{k}_{m}$ with
$$ - \Delta_g Y^{k}_{m} = k(k+1) Y^{k}_{m}, \quad h D_{\theta} Y^{k}_{m} = m Y^{k}_{m}.$$
These eigenfunctions can be separated into various sequences (i.e. {\em ladders} ) associated with different values $ ( \in {\mathcal P}(T^*S^2)$; specifically, the correspondence is given by $c =  \lim_{m \to \infty} \frac{m}{k}).$ The eigenfunctions with maximal $L^\infty$ blow-up are the  sequence of {\em zonal } harmonics given by
\begin{equation} \label{zonal}
u_{h}(r,\theta) = Y^{k}_{0}(r,\theta) = \frac{\sqrt{2k+1}}{2\pi} \int_{0}^{2\pi} ( \cos r + i \sin r \cos \tau )^{k} d\tau;  \quad h = k^{-1}, \, k=1,2,3,...  \end{equation}
It is obvious from (\ref{zonal}) that
$$ |Y^{k}_{0}(p)| \approx k^{1/2}$$
and thus attains the maximal sup growth at $p$ (similarily, at the south pole). 
At the classical level, the zonals $u_h = Y^{k}_{0}$ concentrate microlocally on the Lagrangian tori $\Lambda_{0} = {\mathcal P}^{-1}(1,0)$. From the formula (\ref{polar}) it is clear that away from the poles (where $(r,\theta)$ are honest coordinates),
\begin{equation} \label{reg}
\Lambda_{0} \setminus \{\pm p \}= \{ (r,\theta, \xi_{r} = \pm 1, \xi_{\theta} = 0), \,  r \in (0,\pi) \} \cong S^{2} \setminus \{ \pm p \}. \end{equation}
The choice of $\xi_{r} = \pm 1$ determines the Lagrangian torus (there are two of them) and also, either torus clearly covers the entire sphere. 
At the poles themselves, the projection $\pi_{\Lambda_0}: \Lambda_{0} \to S^2$ has a blowdown singularity with
\begin{equation} \label{blowdown}
 \pi_{\Lambda_0}^{-1}(\pm p)  = S_{\pm}^{*}(S^2) \cong S^1. \end{equation}
To see this, consider the behaviour at $p$ (with a similar computation at $-p$). Rewriting the integral in involution in Euclidean coordinates $(x,y,z) \in \R^3$ one has
$H = (x\xi_y - y \xi_x)^2 +  (x\xi_z - z \xi_x)^2 + (y\xi_z - z \xi_y)^2$
and $\xi_\theta = x \xi_y - y \xi_x.$ Setting $H=1$, $x \xi_y - y \xi_x = 0$ and $(x,y,z) = (0,0,1)$ gives
$$ \pi_{\Lambda_0}^{-1}(p) \cong \{ (\xi_x,\xi_y) \in \R^2; \xi_x^2 + \xi_y^2 = 1 \}.$$
It is then clear from (\ref{reg}) and (\ref{blowdown}) that $\pi_{\Lambda_0} : \Lambda_0 \to S^2$ is surjective and a  diffeomorphism away from the poles (modulo choice of Lagrangian cover) and the fibres above the poles are $S_{\pm}^*(S^2) \cong S^1.$ \red{ We also note that the Lagrangian $\Lambda_0 = \Lambda_{p,2\pi}$ is the $2\pi$-flowout Lagrangian of $S_p^*(S^2)$ and the cylinder $A_p(\delta/2,3\delta)$ is just a local slice of this Lagrangian.}

The defect measure $\mu$ associated with the zonals is

$$ d\mu = |d\theta_1 d\theta_2|,$$
where $(\theta_1,\theta_2; I_1,I_2) \in \R^2/\Z^2 \times \R^2$ are symplectic action-angle variables defined in a neighbourhood of the Lagrangian torus $\Lambda_{0}$ \cite{TZ03}. 
One can choose one of the angle variables $\theta_1 \in S_{p}^*(S^2)$ to parametrize the circle fibre above $p$ (a homology generator of the torus). Then, by the Liouville-Arnold Theorem, the geodesic flow on the torus $\Lambda_0 = \{ I_1 = c_1, I_2 = c_2\}$ is affine with
$$ \theta_j(t) = \theta_j(0) + \alpha_j t, \quad \alpha_j = \frac{\partial H}{\partial I_j} \neq 0.$$
It is then clear that
$$\mu (\Lambda_{p,\delta}) = \int_{0}^{2\pi} d\theta_1 \cdot \int_{|t| < \delta} \alpha_2 dt \approx \delta \neq 0$$
and $\supp \mu|_{\Lambda_p}=\Lambda_p$. Therefore, this case violates the assumption in Theorem~\ref{mainthm} and that is of course consistent with the maximal $L^\infty$ growth of zonal harmonics.

The analysis above extends in a straightforward fashion to the case of a more general sphere of rotation \cite{TZ03}.



\section{Eigenfunctions of Schr\"{o}dinger operators}
Consider a Schr\"{o}dinger operator  $P(h) = -h^{2} \Delta_g + V$ with $V \in C^{\infty}(M;\R)$ on a compact, closed Riemannian manifold $(M,g)$ and let $u_h$ be $L^2$-normalized eigenfunction with
\begin{equation}
\label{e:schrod} P(h) u_h = E(h) u_h, \quad E(h) = E + o(1), \,\, E> \min V,\qquad \|u_h\|_{L^2}=1.
\end{equation}
Any sequence $u_h$ of solutions to~\eqref{e:schrod} has a subsequence $u_{h_k}$ with a defect measure $\mu$ in the sense that for $a\in C_0^\infty(T^*M)$
$$\langle a(x,hD)u_h,u_h\rangle\to \int_{T^*M} ad\mu.$$
Such a measure $\mu$ is supported on $\{p=0\}$ and is invariant under the bicharacteristic flow  $G_t:=\exp(tH_p).$ 

In analogy with the homogeneous case, we define for $x\in M$  the \emph{time $T$ flow out} by
$$\Lambda_{x,T,V}:=\bigcup_{t=-T}^TG_t(\Sigma_x)$$
where 
$$\Sigma_x=\{\xi\in T^*_xM\mid |\xi|_g^2 + V(x)=E\}.$$ 
\begin{defn} 
Let ${\mathcal H}^{n}$ be $n$-dimensional Hausdorff measure on $\{|\xi|_g^2 +V(x)=E\}$ induced by the Sasaki metric on $T^*M$.
 We say that the sequence $u_{h}$ of solutions to~\eqref{e:schrod} is \emph{\admissible at $x$} if \red{for any defect measure $\mu$ associated to the sequence,} \red{there exists $T>0$ so that}
\begin{equation} \label{admissibleV}
{\mathcal H}^n (  \,\supp \,\, \mu |_{\Lambda_{x,\red{T,}V}} \, ) = 0. \end{equation}
\end{defn}
With these definitions we have the analog of Theorem~\ref{thmweakscarring}
\begin{theo} \label{schrodinger}
Let $B \subset V^{-1}(E)$ be a closed ball in the classically allowable region and $\mu$ be a defect measure associated with the eigenfunction sequence $u_h.$  Then, if the eigenfunction sequence is admissible  for all $x\in B$ in the sense of (\ref{admissibleV}),
$$ \sup_{x \in B} |u_h(x)| = o(h^{\frac{1-n}{2}}).$$
 \end{theo}
\begin{proof}
In analogy with the homogeneous case \cite[Lemma 5.1]{CHT}, we have

$$ \rho( h^{-1} [ P(h) -E]))(x,y) =  h^{\frac{1-n}{2}} a(x,y,h) e^{-i A(x,y)/h} + R(x,y,h)$$
where $A(x,y) \in [ (2C_0)^{-1} \epsilon, 2 C_0 \epsilon]$ for some $C_0>1$ and is the action function defined to be the integral of the Lagrangian $L(x,\xi) = |\xi|_g^2 - V(x)$ along the bicharacteristic in $\{ p = E \}$ starting at $(y,\eta)$ and ending at $(x,\xi).$ For $(x,y)$ in a small neighborhood of the diagonal, there is a unique such $\eta$ satisfying this condition.  The remainder $R(x,y,h) = O(h^{\infty})$ pointwise and with all derivatives. The proof then follows using the same argument  as in the homogeneous case. \end{proof}
\bibliography{biblio}
\bibliographystyle{alpha}

\end{document}